\documentclass[12pt]{article}

\usepackage{amssymb,amsmath,amsthm}
\usepackage{amsfonts}
\usepackage{enumerate}

\usepackage{fancyhdr}	
\usepackage{mdframed}

\usepackage{geometry}
\newcommand\blfootnote[1]{%
  \begingroup
  \renewcommand\thefootnote{}\footnote{#1}%
  \addtocounter{footnote}{-1}%
  \endgroup
}

\newtheoremstyle{mytheoremstyle} 
    {10pt}                    
    {10pt}                    
    {\normalfont}                   
    {}                           
    {\bfseries}                   
    {.}                          
    {0.3cm}                       
    {}  
\theoremstyle{mytheoremstyle}

\theoremstyle{plain}

\fancypagestyle{plain}{%
  \fancyhf{}%
}

\newtheorem{theorem}{Theorem}[section]

\newtheorem{lemma}[theorem]{Lemma}
\newtheorem{proposition}[theorem]{Proposition}

\newtheorem{remark}[theorem]{Remark}

\newtheorem{corollary}[theorem]{Corollary}
\newtheorem{problem}[theorem]{Problem}

\begin{document}

\title{Characterization of $k$-spectrally monomorphic Hermitian matrices}
\author{
	Kawtar Attas, Abderrahim Boussaïri, Imane Souktani
}

\maketitle

\begin{abstract}
This paper solves the following problem about Hermitian matrices related to
the theory of $2$-structures:\emph{ }Let $n$ be a positive integer and $k$ be
an integer with $k\in \{3,\ldots,n-3\}$. Characterize the Hermitian matrices
$A$ such that the characteristic polynomials of the $k\times k$ submatrices of
$A$ are all equal. Such matrices are called $k$-spectrally monomorphic.
A crucial step to obtain this characterization is proving that if a matrix $A$ is $k$-spectrally monomorphic
then it is $l$-spectrally monomorphic for $l$ in $\{1,\ldots,min\{k, n-k\}\}$.
\end{abstract}



\section{Introduction}
\blfootnote{Corresponding author: Abderrahim Boussaïri. Email: aboussairi@hotmail.com}
\blfootnote{Laboratoire Topologie, Alg\`ebre, G\'eom\'etrie et Math\'ematiques Discr\`etes, Facult\'e des Sciences A\"in Chock, Hassan II University of Casablanca, Maroc.}
Let $n$ be a positive integer. An $n\times n$ matrix $A$ is $k$%
-\emph{spectrally monomorphic} if all its $k\times k$ submatrices have the
same characteristic polynomial. Spectral monomorphy is closely related to the
notion of monomorphy introduced by Fra\"{\i}ss\'{e} \cite{Fraissi86}. We will
define this notion for \emph{labeled }$2$\emph{-structure}. Following
\cite{Ehren99}, a \emph{labeled }$2$\emph{-structure }on a set $V$, or shortly
an $l2$-structure, is a map $g$ from the set $\{(x,y):x\neq y\in V\}$ to a
\emph{label set} $\mathcal{C}$. The elements of $V$ are called the
\emph{vertices} of $g$. With each subset $X$ of $V$, we associate the
$l2$\emph{-substructure} $g[X]$ of $g$, \emph{induced} by $X$, defined by
$g[X](x,y):=g(x,y)$ for any $x\neq y\in X$. $l2$-structures were introduced to
generalize the notion of graphs, tournaments and other binary structures.
Recall that an $n$\emph{-tournament} $T$ is a digraph with $n$ vertices in
which every pair of distinct vertices is joined by exactly one arc. If the arc
joining vertices $u$ and $v$ of $T$ is directed from $u$ to $v$, then $u$ is
said to \emph{dominate} $v$ (symbolically $u\rightarrow v$). For more details
about tournaments, we refer the reader to \cite{Moon68}.

Let $g$ and $h$ be $l2$-structures with the same label set and whose vertex
sets are, respectively, $V$ and $W$. We say that $g$ and $h$ are
\emph{isomorphic} if there exists a bijection $\sigma$ from $V$ onto $W$ such
that $g(x,y)=h(\sigma(x),\sigma(y))$ for any $x\neq y\in V$. An $l2$-structure
is $k$-\emph{monomorphic} if all its substructures on $k$ vertices are isomorphic.

Several results about monomorphic relations were obtained by Assous
\cite{Roland86}, Frasnay \cite{Frasnay65}, and Pouzet \cite{Pouzet76,
Pouzet81}. A basic example is the class of transitive tournaments.
A tournament $T$ is \emph{transitive} if, whenever $x\rightarrow y$ and
$y\rightarrow z$, then $x\rightarrow z$. The smallest non-transitive
tournament consists of $3$ vertices $x$, $y$ and $z$ such that $x\rightarrow
y\rightarrow z\rightarrow x$. Such a tournament is called a $3$\emph{-cycle}.
It is not difficult to see that $3$-monomorphic tournaments with at least $4$
vertices are transitive. Moreover, it follows from a combinatorial lemma of
Pouzet \cite{Pouzet76} that $k$-monomorphic $n$-tournaments, for
$k=3,\ldots,n-3$, are $3$-monomorphic, and hence are transitive. For $k=n-2$,
an $n$-tournament whose automorphism group acts transitively on the set of its
arcs is $(n-2)$-monomorphic. These tournaments are called \emph{arc-symmetric}
and were characterized independently by Kantor \cite{kantor69} and Berggren
\cite{Berggren72}. Conversely, Jean \cite{Jean69} proved earlier that
$(n-2)$-monomorphic $n$-tournaments with at least $5$ vertices are either
transitive or arc-symmetric. The problem of the characterization of
$(n-1)$-monomorphic $n$-tournaments proposed by Kotzig (see \cite{Bondy76},
problem 43, p. 252) remains unsolved. Some progress was
made by Yucai et al. \cite{Yucai87} and Issawi \cite{Issawi18}. In
\cite{Boudabbous99}, Boudabbous proposed a weak notion of monomorphy for
tournaments using \textquotedblright isomporphy up to
complementation\textquotedblright \ instead of \textquotedblright
isomporphy\textquotedblright. An analogous study for graphs was done by
Boushabi and Boussa\"{\i}ri \cite{Boussairi12}.

Let $g$ be a \emph{complex }$l2$\emph{-structure} with $n$ vertices, that is,
its label set is the complex field. With respect to an ordering $x_{1}%
,\ldots,x_{n}$ of the vertex set, we can identify $g$ to the $n\times n$ zero
diagonal matrix $M=[m_{ij}]$ in which $m_{ij}=g(x_{i},x_{j})$ if $i\neq j$.
This defines a one-to-one correspondence between $n\times n$ zero-diagonal
complex matrices, and $l2$-structures with vertex set $\{x_{1}
,\ldots,x_{n}\}$. Note that a complex $l2$-structure is $k$-monomorphic if
all the $k\times k$ submatrices of its corresponding matrix are
permutationally similar. Hence, $k$-spectral monomorphy for zero-diagonal
complex matrices is a weakening of $k$-monomorphy for $l2$-structures.

In this paper, we address the following problem.

\begin{problem}
\label{original-pb} Characterize $k$-spectrally monomorphic $n\times n$
matrices, where $n$ is a positive integer and $k\in \{3,\ldots,n-1\}$.
\end{problem}

We restrict ourselves to Hermitian matrices. The class of $l2$-structures
corresponding to Hermitian matrices\ encompass many well-known classes\ like
graphs, tournaments, and digraphs.

\section{Some properties of $k$-spectral monomorphy}
We will give some operations that preserve hermiticity and $k$-spectral monomorphy. Let $H$ be
an $n\times n$ Hermitian matrix. The characteristic polynomial of $H$ is
$\phi_{H}(x)=\det(xI-H)$. For a subset $\alpha$ of $\{1,\ldots,n\}$, we denote
by $H[\alpha]$ the principal submatrix of $H$ whose rows and
columns are indexed by $\alpha$.

\begin{itemize}
\item Let $a$ be a real number. Then
\begin{equation}
\phi_{( H-a I)  [  \alpha ]  }(x)=\phi_{H[ \alpha ] }(x+a) \label{eq s}
\end{equation}

\item Let $P=(p_{ij})$ be an $n\times n$ permutation matrix and let $\sigma$
be the corresponding permutation, that is $p_{ij}=1$ if and only if
$\sigma(i)=j$. For every subset $\alpha$ of $\{1,\ldots,n\}$, we have
$PHP^{t}[  \alpha ] =H[  \sigma^{-1}(  \alpha )
]  $, and hence
\begin{equation}
\phi_{PHP^{t}[  \alpha ]  }(x)=\phi_{H[  \sigma^{-1}
(\alpha)]  }(x) \label{Eqb}
\end{equation}

\item Let $D$ be a diagonal matrix whose diagonal entries have modulus $1$.
Then, for every subset $\alpha$ of $\{1,\ldots,n\}$, we have
\begin{equation}
\phi_{DH\overline{D}[  \alpha ]  }(x)=\phi_{H[  \alpha ]
}(x) \label{Eq a}
\end{equation}

\item Let $b$ be a non-zero real number. Then
\begin{equation}
\phi_{b H[  \alpha ]  }(x)=b^{\vert \alpha
 \vert }\phi_{H[\alpha]  }(\tfrac{x}{b}) \label{Eqc}
\end{equation}

\end{itemize}

Let $\Gamma$ be the subgroup of the general linear group, generated by the
unitary diagonal matrices and the permutation matrices.
Two Hermitian matrices $H_{1}$ and $H_{2}$ are $\Gamma$-equivalent
if $H_{2}= a(  SH_{1}S^{*})  + bI$ for some real
numbers $a$ and $b$, and $S$ in $\Gamma$. This defines an equivalence relation
between $n\times n$ Hermitian matrices. It follows from the above equalities
that this relation preserves $k$-spectral monomorphy.

Problem \ref{original-pb} is trivial for $k=1$ and for $k=2$. Indeed, if $H$
is a $n\times n$ Hermitian matrix with $n\geq3$, then

\begin{itemize}
\item[1)] $H$ is $1$-spectrally monomorphic if and only if all its
diagonal entries are equal.

\item[2)] $H$ is $2$-spectrally monomorphic if and only if it is
$1$-spectrally monomorphic and all of its off-diagonal entries have the same modulus.
\end{itemize}

We say that a matrix $A=(  a_{ij})  $ is
\emph{normalized} if $a_{i1}=a_{1i}=1$ for $i\neq1$.

\begin{remark}
\label{remark-module}Let $H$ be an $n\times n$ Hermitian matrix with a
non-zero off-diagonal entry. If $H$ is $2$-spectrally monomorphic, then $H$ is
$\Gamma$-equivalent to a normalized zero-diagonal Hermitian matrix.
\end{remark}

A fundamental property of $k$-spectral monomorphy is given in the following proposition.

\begin{proposition}
\label{prop inf kspectra}If $A$ is\textit{ a }$k$-spectrally monomorphic
\textit{complex matrix, then }$A$\textit{ is }$l$-spectrally monomorphic for
each $l\in \left \{  1,\ldots,\min(k,n-k)\right \}  $.
\end{proposition}

To prove this proposition, we will apply the following result which is a
consequence of \cite[Lemma~II-2.2]{Pouzet76}.

\begin{lemma}
\label{lemm-these20}Let $V$ be a set of size $n$. Let $p$ and $r$ be two
arbitrary integers satisfying $n\geq p+r$ and let $f$ be a map from $\binom
{V}{p}$ to the complex field $\mathbb{\mathbb{C}}$.
If $\underset{P\subseteq B}{\sum}f(P)$ is independent of $B$, where
$B\in \binom{V}{p+r}$, then for any subset $X$ of $V$ such that $|X|\leq
n-(p+r)$, the number $\underset{P\supseteq X}{\sum}f(P)$ depends only on the
cardinality of $X$. Moreover,$\ $if $n\geq2p+r$, then $f$ is a constant map.
\end{lemma}

We need also the following lemma (see for example \cite[p.294]{Meyer}).

\begin{lemma}
\label{lecoeff-polynom} Let $A$ be an $n\times n$ matrix with
characteristic polynomial $\phi(x):=x^{n}+a_{1}x^{n-1}+a_{2}x^{n-2}
+\cdots+a_{n-1}x+a_{n}$. Then
\begin{equation}
a_{p}=(-1)^{p}\underset{\alpha \in \binom{\{1,\ldots,n\}}{p}}{\sum}\det A[\alpha]  \label{minor}
\end{equation}

for $p=1,\ldots,n$.
\end{lemma}

\begin{proof}[Proof of Proposition \ref{prop inf kspectra}]
By Lemma \ref{lecoeff-polynom}, it suffices to prove that for every pair of
subsets $\alpha$ and $\beta$ of $\left \{
1,\ldots,n\right \}  $ such that $|\alpha|=|\beta|\leq \min(k,n-k)$, we have
$\det A[\alpha]=\det A[\beta]$. Let $p\leq \min(k,n-k)$ and let $r:=k-p$. We
will apply Lemma \ref{lemm-these20} to the map $f(\theta):=\det A[\theta]$,
where $\theta \in \binom{\left \{  1,\ldots,n\right \}  }{p}$. For this, let
$\gamma_{1},\gamma_{2}\subseteq \left \{  1,\ldots,n\right \}  $ such that
$|\gamma_{1}|=|\gamma_{2}|=p+r=k$. Since $A$ \textit{is }$k$-spectrally
monomorphic, $A[\gamma_{1}]$ and $A[\gamma_{2}]$ have the same characteristic
polynomial $x^{k}+a_{1}x^{k-1}+\cdots+a_{p}x^{k-p}+\cdots+a_{k}$. From Lemma
\ref{lecoeff-polynom}, we have
\[
a_{p}=(-1)^{p}\underset{\theta \in \binom{\gamma_{1}}{p}}{\sum}\det
A[\theta]=(-1)^{p}\underset{\theta \in \binom{\gamma_{2}}{p}}{\sum}\det
A[\theta]\text{.}
\]
Or, equivalently,
\[
\underset{\theta \in \binom{\gamma_{1}}{p}}{\sum}f(\theta)=\underset{\theta
\in \binom{\gamma_{2}}{p}}{\sum}f(\theta)
\]
Moreover, $n\geq2p+r$ because $p\leq \min(k,n-k)$. Thus, by Lemma
\ref{lemm-these20}, $f$ is a constant map.
\end{proof}

The next corollary is a particular case of Proposition \ref{prop inf kspectra}.

\begin{corollary}
\label{prop-consec}Let $A$ be a $k$-spectrally monomorphic
$n\times n$ complex matrix. If $n\geq2k-1$, then the
following equivalent assertions hold

\begin{description}
\item[1)] For every pair of subsets $\alpha$ and $\beta$ of $\left \{  1,\ldots
,n\right \}  $ with $|\alpha|=|\beta|\leq k$, we have $\det A[\alpha]=\det
A[\beta]$.

\item[2)] The matrix $A$ is $l$-spectrally monomorphic for $l=1,\ldots,k$.
\end{description}
\end{corollary}

\section{\sloppy $k$-spectrally monomorphic Hermitian
matrices for $k\in \{3,\ldots,n-4\}$}

\ \ Let $c$ be a complex number. The Hermitian matrix
\[
H_{n}(  c)  =
\begin{pmatrix}
0 & c & \cdots & \cdots & c\\
\overline{c} & 0 & c & \cdots & c\\
\vdots & \overline{c} & \ddots & \ddots & \vdots \\
\vdots & \vdots & \ddots & 0 & c\\
\overline{c} & \overline{c} & \cdots & \overline{c} & 0
\end{pmatrix}
\]
is $k$-spectrally monomorphic for every $k\in \{1,\ldots,n\}$.

Let $S$ be a skew-symmetric matrix whose off-diagonal entries are from the set
$\{-1,1\}$. The Hermitian matrix $iS$ is $k$-spectrally monomorphic for
$k\in \{1,2,3\}$.

Let $n\geq4$ and let $Q=(  q_{ij})$ be a $3$-spectrally monomorphic $n\times n$
zero-diagonal Hermitian matrix. Assume that there exists a non-real complex
number $c$ of modulus $1$ such that for every $i\neq j\in \{1,\ldots,n\}$,
$q_{ij}\in \{c,\overline{c}\}$ and $q_{1j}=c$ for $j\in \{2,\ldots,n\}$.

\begin{lemma}
Consider the tournament $T$ with vertex set $\{1,\ldots,n\}$ such that $i$
dominates $j$ if and only if $q_{ij}=c$.
\label{lemma-permutation}
\begin{itemize}
\item[i)] If $T$ is transitive, then there exists a permutation matrix $P$
such that $H_{n}(  c)  =PQP^{t}$.

\item[ii)] If $T$ is not transitive, then
$c\in \{i,-i\}$, or equivalently $Q=iS$ for some normalized skew-symmetric
matrix $S$ whose off-diagonal entries are from the set $\{-1,1\}$. Moreover,
$Q$ is not $4$-spectrally monomorphic.
\end{itemize}
\end{lemma}

\begin{proof}
Firstly, assume that $T$ is transitive. There exists a permutation $\sigma$
such that if $i<j$, then $\sigma (  i)  \rightarrow \sigma (
j)  $. It is easy to see that $H_{n}(  c)  =PQP^{t}$, where
$P$ is the permutation matrix corresponding to $\sigma$.

To prove ii), assume that $T$ contains a $3$-cycle $j\rightarrow k\rightarrow
l\rightarrow j$. Then
\[
\det Q[j,k,l]=c^{3}+\overline{c}^{3}
\]
Since $q_{1s}=c$ for $s\in \{2,\ldots,n\}$, $1\notin \{j,k,l\}$. Moreover
\[
\det Q[1,j,k]=c+\overline{c}
\]
As $Q$ is $3$-spectrally monomorphic, $c^{3}+\overline{c}
^{3}=c+\overline{c}$. This implies that $c\in \left \{  i,-i\right \}  $ because
$c\overline{c}=1$ and $c$ is not real. Now, we will prove that $Q$ is not
$4$-spectrally monomorphic. We have $\det Q[1,j,k,l]=9$. However, for
$m\notin \left \{  1,j,k,l\right \}  $, it is not difficult to check that if
$m\rightarrow k$, then $\det Q[1,j,k,m]=1$ and if $k\rightarrow m$, then $\det
Q[1,k,l,m]=1$.
\end{proof}

For $k\in \{3,4\}$, we have the following characterization.

\begin{theorem}
\label{main th-one}Let $H$ be an $n\times n$ Hermitian matrix. Then the
following hold

\begin{itemize}
\item[i)] If $n\geq5$, then $H$ is $3$-spectrally monomorphic if and only if
$H$ is $\Gamma$-equivalent to $H_{n}(  c)  $ where $c$ is a complex
number, or $\Gamma$-equivalent to $iS$ for some\ skew-symmetric matrix $S$
whose off-diagonal entries are from the set $\{-1,1\}$;

\item[ii)] If $n\geq7$, then $H$ is $4$-spectrally monomorphic if and only if
$H$ is $\Gamma$-equivalent to $H_{n}(  c)  $ where $c$ is a complex number.
\end{itemize}
\end{theorem}

To prove this theorem, we start with the study of the possible entries of
normalized $3$-spectrally monomorphic Hermitian matrices.

\begin{lemma}
\label{lem-normalized-u3monom-}Let $H=(  h_{ij})  $ be an\ $n\times
n$ zero-diagonal Hermitian normalized matrix with $n\geq5$. If $H$ is
$3$-spectrally monomorphic, then there exists a complex number $c$ of modulus
$1$ such that for every $i\neq j\in \{2,\ldots,n\}$, we have $h_{ij}
\in \{c,\overline{c}\}$.
\end{lemma}

\begin{proof}
By Corollary \ref{prop-consec}, $\det H[\alpha]=\det H[\beta]$ for every
subsets $\alpha$ and $\beta$ of $\{1,\ldots,n\}$ such that $|\alpha
|=|\beta|\leq3$. Let $i\neq j\in \{2,\ldots,n\}$. Then,
\begin{align*}
\det H[  i,j]   &  =\det H[  1,2] \\
\det H[1,i,j]  &  =\det H[1,2,3].
\end{align*}

It follows that
\[
\left \vert h_{ij}\right \vert =\left \vert h_{12}\right \vert =1
\]
Moreover
\[
\det H[1,i,j]=h_{ij}+\overline{h_{ij}}\mbox{.}
\]

	Let $c:=h_{23}$. Then $h_{ij}+\overline{h_{ij}} =  c+\overline{c}$ and
$h_{ij}\overline{h_{ij}} =  c\overline{c}$. Hence $h_{ij}\in \{c,\overline{c}\}$.
\end{proof}

As a consequence, we obtain the following.
\begin{corollary}
\label{corol-descr2struc}Let $H$ be a $3$-spectrally monomorphic $n\times n$
Hermitian matrix with $n\geq5$. Then $H$ is a real scalar matrix or $H$ is
$\Gamma$-equivalent to a Hermitian zero-diagonal matrix $Q=(
q_{ij})  $ such that $q_{ij}\in \{c,\overline{c}\}$ for $i\neq j$ and
$q_{1j}=c$ for $j\in \{2,\ldots,n\}$, where $c$ is a complex number of modulus
$1$.
\end{corollary}

\begin{proof}
By the second assertion of Corollary \ref{prop-consec}, $H$ is $k$-spectrally
monomorphic for $k=1$ and $k=2$. If $H$ is a diagonal matrix, then $H$ is a
real scalar matrix. Assume that $H$ has a non-zero off-diagonal entry. It
follows from Remark \ref{remark-module} that $H$ is $\Gamma$-equivalent to a
normalized zero-diagonal Hermitian matrix $\widetilde{H}=(  \widetilde
{h}_{ij})  $. As the $\Gamma$-equivalence preserves the spectral
monomorphy, the matrix $\widetilde{H}$ is $3$-spectrally monomorphic. Thus, by
Lemma \ref{lem-normalized-u3monom-}, there exists a complex number $c$ of
modulus $1$ such that for every $i\neq j\in \{2,\ldots,n\}$, we have
$\widetilde{h}_{ij}\in \{c,\overline{c}\}$. To conclude, it suffices to choose
$Q=\overline{D}\widetilde{H}D$, where $D=diag(  1,c,\ldots
,c)  $.
\end{proof}

Now, we are able to prove Theorem \ref{main th-one}.

\begin{proof}[Proof of Theorem \ref{main th-one}]
i) If $H$ is a real scalar matrix, then $H$ is
$\Gamma$-equivalent to $H_{n}(  0)  $. Assume that $H$ is not a
real scalar matrix, and consider the matrix $Q$, $\Gamma$-equivalent to $H$,
as described in Corollary \ref{corol-descr2struc}. If $c\ $is a real number,
then $c\in \left \{  1,-1\right \}  $, and hence $Q=\pm H_{n}(  1)  $.
If $c\ $is not a real number, then the result is obtained by applying Lemma
\ref{lemma-permutation}.

ii) By Corollary \ref{prop-consec}, $H$ is $3$-spectrally monomorphic. Using
assertion i) above, we can assume that $H$ is $\Gamma$-equivalent to $iS$ for
some normalized skew-symmetric matrix $S$ whose off-diagonal entries are from
the set $\{-1,1\}$. By Lemma \ref{lemma-permutation}, the matrix $iS$ is
permutationally similar to $H_{n}(  i)  $. Hence $H$ is $\Gamma
$-equivalent to $H_{n}(  i)  $.
\end{proof}

The following corollary is a direct consequence of Proposition
\ref{prop inf kspectra} and assertion ii) of Theorem \ref{main th-one}.

\begin{corollary}
Let $H$ be an $n\times n$ Hermitian matrix. If $n\geq8$ and $4\leq k\leq n-4$,
then $H$ is $k$-spectrally monomorphic if and only if $H$ is $\Gamma
$-equivalent to $H_{n}(  c)  $ where $c$ is a complex number.
\end{corollary}

\begin{remark}
The classes of $k$-spectrally monomorphic Hermitian matrices for $k=3$ and
$n=4$, as well as $k=4$ and $n=5,6$, are not easy to describe. More generally,
the characterization of $k$-spectrally monomorphic $n\times n$ Hermitian
matrices for $k=n-1,n-2$ seems difficult.
\end{remark}

\section{Characterization of $(n-3)$-spectrally monomorphic Hermitian
matrices}

\ \ As we have seen above, the Hermitian matrices $H_{n}(  c)  $ are
$(n-3)$-spectrally monomorphic. Another example is obtained from
skew-symmetric conference matrices. Recall that a \emph{conference matrix} is
an $n\times n$ matrix $C$ with $0$ on the diagonal and $1$ and $-1$ off the
diagonal, such that $C^{t}C=(  n-1)  I_{n}$. Let $S$ be a skew-symmetric conference matrix.
Assertion iv) of the next proposition shows that the Hermitian matrix $iS$ is
$(n-3)$-spectrally monomorphic.

\begin{proposition}
\label{hermitian spectral}Let $S$ be a skew-symmetric conference matrix of
order $4t+4$. Then

\begin{itemize}
\item[i)] The characteristic polynomial of $iS$\ is $(x^{2}-4t-3)^{2t+2}$.

\item[ii)] The characteristic polynomial of the matrix obtained from $iS$\ by
deleting one row and the corresponding column is$\ x(x^{2}-4t-3)^{2t+1}$.

\item[iii)] The characteristic polynomial of the matrix obtained from $iS$\ by
deleting two rows and the corresponding columns is $(x^{2}-1)(x^{2}%
-4t-3)^{2t}$.

\item[iv)] The characteristic polynomial of the matrix obtained from $iS$\ by
deleting three rows and the corresponding columns is $x(x^{2}-3)(x^{2}%
-4t-3)^{2t-1}$.
\end{itemize}
\end{proposition}

The proof of this proposition is contained implicitly in \cite{Suda17}. It is
based on the interlacing theorem due to Cauchy \cite{fisk05}.

The characterization of $(n-3)$-spectrally monomorphic Hermitian matrices is given by
the following theorem.

\begin{theorem}
\label{main th-deux}
Let $H$ be an $n\times n$ Hermitian matrix. If $n\geq7$, then $H$ is
$(n-3)$-spectrally monomorphic if and only if $H$ is $\Gamma$-equivalent to
$H_{n}(  c)  $ where $c$ is a complex number or to $iS$ where $S$
is a skew-symmetric conference matrix.
\end{theorem}

Before proving this theorem, we recall some properties of skew-symmetric
conference matrices and their relationship with doubly regular tournaments.

Skew-symmetric conference matrices are related to skew Hadamard matrices. A
\emph{Hadamard matrix} $H$ is a square matrix of order $n$ whose entries are
from $\left \{  -1,1\right \}  $ and whose rows are mutually orthogonal, or
equivalently, $HH^{t}=H^{t}H=nI_{n}$. The order of a Hadamard matrix is
necessarily $1$, $2$ or a multiple of $4$. It is conjectured \cite{Paley33}
that Hadamard matrices of order $n$ always exist when $n$ is divisible by $4$.
A Hadamard matrix $H$ of order $n$ is called \emph{skew} if $H+H^{t}=2I_{n}$.
It is easy to see that $H$ is a skew Hadamard matrix if and only if $H-I_{n}$
is a skew-symmetric conference matrix. Reid and Brown \cite{Reid72} gave a
construction of skew Hadamard matrices from doubly regular
tournaments.\ Recall that a tournament $T$ of order $n$ is \emph{doubly
regular} if there exists $t>0$, such that every pair of vertices is dominated
by exactly $t$ vertices. We have necessarily  that $n=4t+3$. Let
$\widehat{T}$ be the tournament obtained from $T$ by adding a new vertex which
dominates every vertex of $T$. If $A$ is the adjacency matrix of $\widehat{T}%
$, then $A-A^{t}+I_{4t+4}$ is a skew Hadamard matrix and hence $A-A^{t}$ is a
skew-symmetric conference matrix. Conversely, let $H$ be a normalized skew
Hadamard matrix of order $4t+4$ and let $K$ be the matrix obtained from $H$ by
removing the first row and the corresponding column. We denote by $J_{4t+3}$
the all-ones matrix. Reid and Brown \cite{Reid72} showed that the tournament
with adjacency matrix $\frac{1}{2}(K+J_{4t+3}-2I_{4t+3})$ is doubly regular.

Let $T$ be a tournament and let $i,j$ be two vertices of $T$. We denote by
$C_{3}(i,j)$ (resp. $O_{3}(i,j)$), the number of $3$-cycles (resp. the number
of transitive $3$-tournaments) of $T$ containing $i$ and $j$. The tournament
$T$ is \emph{homogeneous} if there exists an integer $k>0$ such that
$C_{3}(i,j)=k$ for every vertices $i,j$ of $T$. Kotzig \cite{Kotzig69}
proved that such a tournament contains exactly $4k-1$ vertices. Moreover, Reid
and Brown \cite{Reid72} established that it is doubly regular.

\begin{proof}[Proof of Theorem \ref{main th-deux}]
It suffices to prove the direct implication.
By Proposition \ref{prop inf kspectra}, $H$ is $3$-spectrally monomorphic.
Using the first assertion of Theorem \ref{main th-one}, we can assume that $H$
is $\Gamma$-equivalent to $iS$ where $S=(  s_{ij})  $ is a
skew-symmetric matrix. Without loss of generality, the matrix $S$ can be
chosen to be normalized. Consider the tournament $T$ with vertex set $\left \{
2,\ldots,n\right \}  $ such that $i$ dominates $j$ if $s_{ij}=1$. We can assume
that $T$ is not transitive, because otherwise $iS$ is $\Gamma$-equivalent to
$H_{n}(  i)  $. We have to prove that $S$ is a skew-symmetric
conference matrix, or equivalently, $T$ is homogeneous.

Consider two arbitrary vertices $i,j$ of $T$. It is easy to see that for
every $k\in \left \{  2,\ldots,n\right \}  \setminus \{i,j\}$, $\det S[1,i,j,k]=9$
if $i,j$ and $k$ form a $3$-cycle of $T$ and $\det S[1,i,j,k]=1$ otherwise.
It follows that%
\[
\underset{\substack{|\alpha|=4 \\ \{1,i,j\} \subseteq \alpha}}{\sum}\det S[\alpha
]=9\cdot C_{3}(i,j)+O_{3}(i,j)
\]

Since $C_{3}(i,j)+O_{3}(i,j) = n-3$, we have
\[
8\cdot C_{3}(i,j)=\underset{\substack{|\alpha|=4\\ \{1,i,j\} \subseteq \alpha}}{\sum
}\det S[\alpha]-(n-3)
\]

To conclude, it suffices to prove that the number 
\[\underset{\substack{|\alpha|=4\\ \{1,i,j\} \subseteq \alpha}}{\sum}\det S[\alpha]\] does not depend on $i$ and
$j$. For this, we will use Lemma \ref{lemm-these20} applied to the map
$f(\alpha):=\det S[\alpha]$ with $V:=\{1,\ldots,n\}$, $p:=4$ and $r:=n-7$. Let $\gamma_{1}$ and
$\gamma_{2}$ be two subsets of $\left \{  1,\ldots,n\right \}  $ with
$|\gamma_{1}|=|\gamma_{2}|=n-3$. Since the matrix $iS$ is $(n-3)$-spectrally
monomorphic, the submatrices $S[\gamma_{1}]$ and $S[\gamma_{2}]$ have the same
characteristic polynomial
\[
\phi(x)=x^{n-3}+a_{1}x^{n-4}+\cdots+a_{p}x^{n-3-p}+\cdots+a_{n-3}
\]
It follows from Lemma \ref{lecoeff-polynom} that
\[
a_{4}=\underset{\substack{|\alpha|=4\\\alpha \subseteq \gamma_{1}}}{\sum}\det S[\alpha
]=\underset{\substack{|\alpha|=4\\ \alpha \subseteq \gamma_{2}}}{\sum}\det S[\alpha]
\]
Hence by Lemma \ref{lemm-these20}, for a subset $\zeta$ of $\left \{
1,\ldots,n\right \}  $ with $|\zeta|=3$, the number $\underset{\substack{|\alpha
|=4\\ \zeta \subseteq \alpha}}{\sum}\det S[\alpha]$ does not depend on $\zeta$.
\end{proof}

\end{document}